\documentclass[11pt]{amsart}
\usepackage{amssymb,amsmath}
\usepackage{amsthm,amscd}
\usepackage{palatino}

\newtheorem{theorem}{Theorem} 
  
\newtheorem{corollary}[theorem]{Corollary} 
\newtheorem{proposition}[theorem]{Proposition}

\theoremstyle{definition}

\newtheorem{remark}[theorem]{Remark}

\newcommand{\g}{\mathfrak{g}}  
\newcommand{\so}{\mathfrak{so}}
\newcommand{\lphi}{\mathfrak{L}\phi}
\newcommand{\tphi}{\tilde{\phi}}
\newcommand{\R}{\mathbf{R}}
\newcommand{\n}{\mathbf{n}}
\newcommand{\ra}{\rightarrow}
\newcommand{\ad}{Ad}

\begin{document}

\title{\textbf{A Connection Whose Curvature is the Lie Bracket}}
\author{ Kent E. Morrison } 
\address{Department of Mathematics \\
         California Polytechnic State University \\
         San Luis Obispo, CA  93407}
 \email { kmorriso@calpoly.edu} 
\date{August 24, 2008}

\begin{abstract}  Let $G$ be a Lie group with Lie algebra $\g$. On the 
trivial principal $G$-bundle over $\g$ there is a natural connection whose 
curvature is the Lie bracket of $\g$. The exponential map of $G$ is given by parallel transport of this connection. If $G$ is the diffeomorphism group of a manifold $M$, the curvature of the natural connection is the Lie bracket of vectorfields on $M$. In the case that $G=SO(3)$ the motion of a sphere rolling on a plane is given by parallel transport of a pullback of the natural connection by a map from the plane to $\so(3)$. The motion of a sphere rolling on an oriented surface in $\R^3$ can be described by a similar connection.
\end{abstract}
\maketitle
\large
\renewcommand{\baselinestretch}{1.2}  \normalsize

\section{A Natural Connection and Its Curvature}

Samelson \cite{Samelson89} has shown that the covariant derivative of 
a connection can be expressed as a Lie bracket. It is the purpose of this 
article to show that the Lie bracket of a Lie algebra can be expressed as the curvature form of a natural connection.

The setting for this result is the following. Let $\pi: P \ra X$ be a right principal $G$-bundle with the Lie group $G$ as the structure group. A connection on $P$ is a smooth $G$-equivariant distribution of horizontal spaces in the tangent bundle $TP$ complementary to the vertical tangent spaces of the fibers. The curvature of a connection is a $\g$-valued 2-form on the total space $P$. A good reference is Bleecker's book \cite{Bleecker81}.

Let $G$ be a Lie group and $\g$ its Lie algebra. Let $P=\g \times G$ be the total space of the trivial right principal $G$-bundle with projection $P \ra \g : (x,g) \mapsto x.$  The right action of $G$ on $P$ is given by $(x,h) \cdot g = (x,hg)$ and let $R_g$ denote the action of $g$ on $P$; that is, $R_g:P \ra P: (x,h) \mapsto (x,hg)$.  Let $1 \in G$ be the identity element and let $\iota: \g \ra \g \times G: x \mapsto (x,1)$ be the identity section of the bundle.

 Let $L_g:G \ra G: h \mapsto gh$ and $R_g:G \ra G : h \mapsto hg$ be the left and right multiplication by $g$. The context should make it clear whether $R_g$ is acting on $P$ or on $G$.
 The adjoint action of $G$ on $\g$ is the derivative at the identity of the conjugation action of $G$ on itself. That is,  
\[ \ad_g :\g \ra \g : v \mapsto T_1 ( L_{g} \circ R_{g^{-1}}) (v) . \]

Finally, we recall the definition of a fundamental vectorfield on a principal $G$-bundle $P$. For $\xi \in \g$, let $\tilde{\xi}$ be the vectorfield on $P$ whose value at $p$ is given by
\[ \tilde{\xi}(p) = \frac{d}{dt} \big(p \cdot \exp(t\xi)\big)\big|_{t=0} .\]
The vectorfield $\tilde{\xi}$ is called a fundamental vectorfield. In the case that $P$ is the trivial bundle $\g \times G$, it is routine to check that 
\[\tilde{\xi}(x,g) = (0,T_1 L_g (\xi)) \in \g \times T_g G. \]

\begin{theorem} There is a natural connection on the trivial bundle $P=\g \times G$ whose local curvature 2-form (with respect to the identity section $\iota$) is the constant $\g$-valued 2-form on $\g$ whose value on a pair of tangent vectors $\xi, \eta \in \g$ is the Lie bracket $[\xi, \eta]$.
\end{theorem}
\begin{proof}
We define the connection on $P$ by choosing the horizontal spaces.
For  $(x,g) \in P$ let $H_{(x,g)} $ be the subspace of  $T_{(x,g)}P$ defined by
      \[
		   H_{(x,g)} := \{ (v,T_1 R_g(v) ) | v \in \g \}.
      \]
It is easy to check that the distribution is right equivariant and that $H_{(x,g)}$ is complementary to the vertical tangent space at $(x,g)$.

Let $\alpha$ be the $\g$-valued 1-form on $P$ defining this connection.  It is characterized by having the horizontal space $H_{(x,g)}$ as the kernel of $\alpha_{(x,g)}$ and by satisfying the conditions
\[ (i)  \quad \alpha_{(x,g)} (\tilde{\xi}(x,g))=\xi , \]
\[(ii) \quad {R_g}^* \alpha = \ad_{g^{-1}} \circ \alpha .\]
 From these properties one can check that $\alpha$ is defined by 
\[
 \alpha_{(x,g)} (v,\xi) = T_g L_{g^{-1}} (\xi) - \ad_{g^{-1}}(v).
\]
The curvature of a connection is the $\g$-valued 2-form on the total space $P$ defined  by 
\[   
   d\alpha + \frac{1}{2}[\alpha, \alpha].
\]
Pulling $\alpha$ back to $\g$ by $\iota$ gives the local connection 1-form $\omega=\iota^*\alpha$. 
 
 The local curvature form $\Omega$ is then
\[ 
  \Omega = d\omega + \frac{1}{2}[\omega , \omega ] .
\]
Computing $\omega$ we see that
\begin{align*}
	\omega_x (v) & =  \alpha_{(x,1)}(v,0)
	\\
	 & =  -v .
\end{align*}
Hence $\omega$ is a constant form, and $d \omega = 0$.
We evaluate $[\omega,\omega]$ on a pair of tangent vectors
$\xi, \eta  \in T_x \g \cong \g$ as follows:
\begin{align*}
   [\omega, \omega ] (\xi,\eta) & =  [\omega (\xi), \omega (\eta) ] - 
   [\omega (\eta), \omega (\xi) ]  \\
    & =  [-\xi,-\eta] - [-\eta,-\xi] \\
    & =  2[\xi,\eta].
\end{align*}
Therefore, the local curvature form $\Omega$ is the constant $\g$-valued 2-form that
maps a pair of tangent vectors $\xi,\eta$, which are just elements of $\g$, to 
the Lie bracket of $\xi$ and $\eta$. 
\[
\Omega(\xi,\eta)= (d \omega  + \frac{1}{2}[\omega,\omega] )(\xi,\eta) = [\xi,\eta]
\]    
\end{proof}

For a Lie group homomorphism $\phi:G_1 \ra G_2$, let $\lphi: \g_1 \ra \g_2$ be the associated Lie algebra homomorphism. (Note that $\lphi$ is $T_1\phi$ when the Lie algebras are viewed as the tangent spaces at the identities of the groups.) Then the map
\[  \tphi: \g_1 \times G_1 \ra \g_2 \times G_2 : (x,g) \mapsto (\lphi(x),\phi(g))\]
is a morphism of principal bundles, which means that it commutes with the right actions of $G_1$ and $G_2$:
\[   \tphi(x,hg)=(\lphi(x),\phi(hg))=(\lphi(x),\phi(h)\phi(g)) .\]
Furthermore, it is a morphism that preserves the horizontal spaces of the natural connections defined in Theorem 1. More precisely, $T\tphi$ maps the horizontal subspaces in $T(\g_1 \times G_1)$ to the horizontal subspaces in $T(\g_2 \times G_2)$ as follows. Let $(v, T_1R_g(v))$ be in $H_{(x,g)} \subset T_{(x,g)}(\g \times T_g G_1)$. 
\begin{align*}  T\tphi(v, T_1R_g(v))  &= (T_x \lphi(v), T_g \phi(T_1 R_g(v))) \\
                                                                 &= (\lphi(v), T_1(\phi \circ R_g)(v)) \\
                                                                 &=(\lphi(v), T_1( R_g \circ \phi)(v)) \\
                                                                 &=(\lphi(v), T_1 R_{\phi(g)}(T_1 \phi(v))),
\end{align*}
which is an element of the horizontal space at $(\lphi(x), \phi(g))$ in $T(\g_2 \times G_2)$.

We also note that for a composition of Lie group homomorphisms $\phi \circ \psi$, the principal bundle map $\widetilde{\phi \circ \psi} = \tphi\circ \tilde{\psi}$, and that $\tilde{I}$ is the identity on the principal bundle for the identity $I: G \ra G$. Hence, we have proved the following theorem.
\begin{theorem}
There is a covariant functor $\mathcal{F}$ from the category of Lie groups to the category of principal bundles with connection, which is defined on objects so that $\mathcal{F}(G)$ is the trivial principal bundle $\g \times G$ with its natural connection and defined on morphisms by $\mathcal{F}(\phi) = \tphi$. 
\end{theorem}

Next we consider the relationship between the connection 1-forms of the two bundles. 
\begin{proposition} The following diagram commutes
\[ \begin{CD}
 T(\g_1 \times G_1)    @> T \tphi >>     T(\g_2 \times G_2) \\
  @V \alpha_1 VV    @VV \alpha_2 V \\
  \g_1 @> \lphi >> \g_2
\end{CD}  \]
Equivalently,
\[
  \tphi^* \alpha_2 = \lphi \circ \alpha_1  .
\]
\end{proposition}
\begin{proof}
Starting with the composition $\alpha_2 \circ T\tphi$ we have
\begin{align*}
(\alpha_2)_{\tphi(x,g)}T\tphi(v,\xi)  &= (\alpha_2)_{(\lphi(x),\phi(g))}(\lphi(v),T_g \phi(\xi)) \\
     &= T_{\phi(g)}L_{\phi(g)^{-1}} (T_g\phi (\xi)) -Ad_{\phi(g)^{-1}}(\lphi(v)) \\
     &= T_g(L_{\phi(g)^{-1}} \circ \phi)(\xi) -Ad_{\phi(g)^{-1}}(\lphi(v)) \\
     &= T_g(\phi \circ L_{g^{-1}} )(\xi) -Ad_{\phi(g)^{-1}}(\lphi(v)) \\
     &= T_1 \phi(T_g L_{g^{-1}}(\xi) - \lphi(Ad_{g^{-1}}(v)) \\
     &= \lphi(T_g L_{g^{-1}}(\xi) - \lphi(Ad_{g^{-1}}(v)) \\
     &= \lphi(T_g L_{g^{-1}}(\xi) -Ad_{g^{-1}}(v) \\
     &= \lphi((\alpha_1)_{(x,g)}(v,\xi) .
\end{align*}
\end{proof}
\begin{proposition} The following diagram commutes
\[ \begin{CD}
 T\g_1     @> T (\lphi) >>     T\g_2 \\
  @V \omega_1 VV    @VV \omega_2 V \\
  \g_1 @> \lphi >> \g_2
\end{CD}  \]
Equivalently,
\[
  (\lphi)^* \omega_2 = \lphi \circ \omega_1.
\]
\end{proposition}
\begin{proof}
Because $\lphi$ is linear, $T_x (\lphi) = \lphi$ for all $x \in \g_1$. Thus, for $v \in T_x \g_1 = \g_1$, 
\[ \omega_2 (T_x(\lphi)(v))=\omega_2(\lphi(v))=-\lphi(v)= \lphi(-v)=\lphi(\omega_1 (v)). \]
\end{proof}
Also, the local curvature forms commute with the Lie algebra homomorphism $\lphi$.
\begin{proposition}
Let $\Omega_i$, $i=1,2$, be the local curvature 2-form for the natural connection on $\g_i \times G_i$. Then $\lphi \circ \Omega_1 = \Omega_2 \circ \lphi.$
\end{proposition}
\begin{proof}
This is simply restating the fact that $\lphi$ is a homomorphism of Lie algebras: $\lphi([\xi,\eta])=[\lphi(\xi),\lphi(\eta)]$.
\end{proof}
\begin{remark}
Let $G$ be the diffeomorphism group of a manifold $M$ with $\g$ being the space of vectorfields on $M$. Although $G$ is not, strictly speaking, a Lie group, the natural connection on $\g \times G$ still makes sense, and so the curvature of this connection is given by the Lie bracket of vectorfields.
\end{remark}

Given a smooth curve $c:[0,1] \ra \g$, \textbf{parallel transport} along $c$ is horizontal lift $(c(t),g(t))$ with initial condition $g(0)=1$. Therefore, $g$ is a solution of the differential equation 
\[ g'(t)= T_1 R_{g(t)} (c'(t)) , \] which simply says that $(c'(t),g'(t))$ is in the horizontal subspace at the point $(c(t), g(t))$. 
  
\begin{theorem} Let $\xi$ be an element of the Lie algebra $\g$ and define $c(t)=t\xi$. Then parallel transport along $c$ is given by $g(t)= \exp(t\xi)$.
\end{theorem}
\begin{proof}
It suffices to show that $g(t)=\exp(t\xi)$ satisfies the differential equation $g'(t)= T_1 R_{g(t)} (c'(t))$ with initial condition $g(0)=1$. The derivative of $\exp(t\xi)$ is given by
\begin{align*} \frac{d}{dt} \exp(t \xi)  &= \frac{d}{ds} \exp((t+s)\xi) |_{s=0} \\
      &= \frac{d}{ds}\left( \exp(s\xi)\exp(t\xi)\right)|_{s=0} \\
     &= \frac{d}{ds}\left(R_{\exp(t\xi)} \exp(s\xi)\right)|_{s=0} \\
      &= T_1R_{\exp(t\xi)}(\xi) \\
      &= T_1 R_{g(t)} (c'(t))
\end{align*}
\end{proof}
\begin{remark} With this result there is another way to see that the Lie bracket is the curvature, since 
\[ [\xi, \eta]
 =\frac{1}{2}\frac{d^2}{dt^2}\exp(t\xi)\exp(t\eta)\exp(-t\xi)\exp(-t\eta)\Big|_{t=0}
 \] 
is the infinitesimal parallel transport around the parallelogram spanned by $\xi$ and $\eta$, which, in turn, is the  curvature tensor applied to $\xi$ and $\eta$.
\end{remark}
\begin{corollary}
Let $\eta$ and $\xi$ be elements of $\g$ and define $c(t)=\eta + t \xi$. Let $g(t)$ be the horizontal lift of $c$ with $g(0)=g_0 \in G$. Then $g(t)=\exp(t\xi) g_0$.
\end{corollary}
\begin{proof}  Since $c'(t)=\xi$, the value of $\eta$ does not matter. Therefore, let $\eta=0$. By the theorem $\gamma(t)=\exp(t\xi)$ is the horizontal lift of $c$ with $\gamma(0)=1$. The right-invariance of the horizontal spaces implies that $g(t) = \gamma(t)g_0$. 
\end{proof}
\begin{remark}
Parallel transport along $c$ is also known as the time-ordered (or path-ordered) exponential of $c'$. One of the equivalent ways to define the time-ordered exponential of a curve $a(t) \in \g$ is  to define it as the solution of the differential equation $g'(t)= T_1 R_{g(t)} (a(t))$ with $g(0)=1$. For $a(t)=c'(t)$ this is the differential equation defining parallel transport along $c$.  To understand the use of the phrase ``time-ordered exponential,'' we use a piecewise linear approximation to $c$ starting at $c(0)$ and consisting of the line segments connecting $c(t_{i-1})$ and $c(t_i)$ where $t_i=(i/n)t$, $i=0,1,\ldots,n$. Let $\Delta t = t/n$. Then by repeated use of the corollary $g(t)$ is approximated by
\begin{align*}
g(t) & \approx \exp(c(t_n)-c(t_{n-1})) \exp(c(t_{n-1})-c(t_{n-2})) \cdots \exp(c(t_1)-c(t_0))  \\
& \approx \exp(\Delta t \,c'(t_{n-1})) \exp(\Delta t \,c'(t_{n-2}))\cdots \exp(\Delta t \,c'(t_0)).
\end{align*}
Therefore, $g(t)$ is the limit as $n$ goes to infinity of this product of exponentials, which are ordered according the parameter value.
\end{remark}

\begin{remark}
The holonomy subgroup (at a point $x_0$ in the base space) of a connection on a principal $G$-bundle is the subgroup of $G$ consisting of the results of parallel transporting around closed curves starting and ending at $x_0$. The null holonomy group is the subgroup resulting from transporting around null-homotopic curves. By the Ambrose-Singer Theorem \cite{AmbroseSinger53} the Lie algebra of the null holonomy group is generated by the values of the curvature tensor. Now $\g$ is simply connected and so the null holonomy group is the holonomy group and its Lie algebra is the derived algebra $[\g, \g]$. Assuming $G$ is connected, the holonomy group is the derived group of $G$.
\end{remark}

\section{Examples}
\subsection{The Special Orthogonal Group}
Let $G$ be the rotation group  $SO(3)$ with  $\g = \so(3)$, the Lie algebra of infinitesimal rotations. We identify $\so(3)$ with $\R^3$ in the standard way so that the vector $v$ in $\R^3$ corresponds to the infinitesimal rotation with axis $v$ that goes counter-clockwise with respect to an observer with $v$ pointing towards him. The angular velocity is the magnitude of $v$. If $v=(v_1,v_2,v_3)$, then the corresponding matrix $\rho_v$ in $\so(3)$ is
\[ \rho_v= 
\left[
\begin{array}{rrr}
 0 &  -v_3 &  v_2 \\
 v_3 & 0  &  -v_1 \\
-v_2  & v_1  &  0 
\end{array}
\right].
\]
Then $\rho_{v \times w} = [\rho_v, \rho_w]$, so that the map $\rho$ is an isomorphism of Lie algebras $(\R^3, \times) \ra (\so(3),[\, , \,])$. 

Given a curve $c:[0,1] \ra \R^3$, parallel transport along $c$ is given by the curve $g:[0,1] \ra SO(3)$, which is the unique solution to the differential equation
\[  g'(t)= \rho_{c'(t)} g(t), \; g(0)=I .\]
In other words, the infinitesimal rotation at $t$ has for its axis of rotation the vector $c'(t)$. One can visualize this as a sphere of radius 1 rotating at the head of a screw (with right hand threads) that is tunneling through space following the trajectory given by the curve $c$. Note that the spherical head is not rigidly attached to the screw because the axis of rotation must be free to vary.

\subsection{A Rolling Sphere on a Plane} 

A variation of this natural connection can be used to describe the geometry of a sphere rolling without slipping on a horizontal plane. The plane is $\R^2$ embedded in $\R^3$ as the set of points $\{(x_1,x_2,0)\}$ and a sphere of radius 1 sits on top of the plane.  Let $c:[0,1] \ra \R^2$ be a smooth curve with initial point $c(0)=(x_1(0),x_2(0))$. Roll the sphere along the curve $c$ until it reaches the endpoint $c(1)=(x_1(1),x_2(1))$. As the sphere rolls along, the point of contact is $c(t)$ and the configuration of the sphere is given by a curve $g(t)$ in $SO(3)$. At each point $c(t)$ the infinitesimal rotation is about the axis $J(c'(t))$, where $J: \R^2 \ra \R^2$ is $\pi/2$ rotation counter-clockwise defined by $J(x_1,x_2)=(x_2,-x_1)$. We can formulate the differential equation satisfied by $g(t)$ as
\[  g'(t)=\rho_{J(c'(t))} g(t), \quad g(0)=I .\]
(Briefly, $g(t + dt) = (I+\rho_{J(c'(t))})\, g(t)\, dt$, from which the differential equation follows.)
In order to see this differential equation as the parallel transport equation for the curve $c$ with initial condition $g(0)=I$, we define a connection on $\R^2 \times SO(3)$ with horizontal space at the point $(x,g)$ in $\R^2 \times SO(3)$ given by
\[H_{(x,g)}=\{(v, \rho_{J(v)}g) | v \in \R^2\} \subset \R^2 \times T_g SO(3). \]
The connection 1-form $\alpha$ is defined by
\[ \alpha_{(x,g)}(v,\xi)= g^{-1} \xi  - g^{-1}\rho_{J(v)} g.\]
Let $\omega = \iota^* \alpha$ be the local connection 1-form. Since
\[ \omega_x(v)=\alpha_{(x,I)}(v,0)=-\rho_{J(v)},\]
we see that $\omega$ is constant, i.e., $d \omega=0$. Then the local curvature 2-form acts on a pair of tangent vectors $u, v \in T_x \R^2=\R^2 \subset \R^3$ by 
\begin{align*} (d\omega + \frac{1}{2}[\omega,\omega])(u, v) &= 0 + \frac{1}{2}[\omega,\omega](u,v) \\
&= \frac{1}{2}( [\omega(u),\omega(v)] - [\omega(v), \omega(u)] ) \\
&= [\omega(u),\omega(v)] \\
&= [-\rho_{J(u)},-\rho_{J(v)}] \\
&= [\rho_{J(u)},\rho_{J(v)}] \\
&=[\rho_u,\rho_v].
\end{align*}
For the last step note that $[\rho_{J(u)},\rho_{J(v)}]=\rho_{J(u)\times J(v) }$ because $\rho$ is an isomorphism of Lie algebras. Also, $J(u)\times J(v) =u \times v $ from the  geometric properties of the cross product . Finally, $\rho_{u \times v} = [\rho_u,\rho_v]$.

The derived algebra of $\so(3)$ is $\so(3)$ because each of the basis elements $e_1, e_2, e_3$ in $\R^3$ is a cross product. The group $SO(3)$ is connected, and so, by the Ambrose-Singer Theorem, the holonomy subgroup is all of $SO(3)$. In other words, any rotation of a sphere can be achieved by rolling the sphere around some closed path in the plane. An elementary proof (without the apparatus of modern differential geometry) of this old result has recently appeared \cite{Johnson07}.

The connection just defined is actually just a pull-back of the natural connection on $\so(3) \times SO(3)$. 
\begin{theorem} Define 
\[ f :\R^2 \ra \so(3): x \mapsto \rho_{J(x)}= \left[
\begin{array}{rrr}
 0 &  0 &  x_1 \\
 0& 0  &  x_2 \\
-x_1  & x_2  &  0 
\end{array}
\right].\]
Then the connection on $\R^2 \times SO(3)$ associated to the rolling sphere is the pullback by $\rho_J$ of the natural connection on $\so(3) \times SO(3)$.
\end{theorem} 
\begin{proof}
The pullback bundle of the trivial bundle is also trivial. Let $\hat{f}$ denote the bundle map 
$\R^2 \times SO(3) \ra \so(3) \times SO(3): (x,g) \mapsto (f(x), g)$. Then it suffices to compute ${\hat{f}}^*\alpha$, the pullback of the connection 1-form $\alpha$ on $\so(3) \times SO(3)$, to see that it is the connection 1-form of the rolling sphere. At a point $(x,g) \in \R^2 \times SO(3)$ we have
\begin{align*}
  (\hat{f}^*\alpha)_{(x,g)} (v,\xi) &= \alpha_{(f(x),g)}(Df(x)(v),\xi) \\
                                                    &= \alpha_{(f(x),g)}(f(v),\xi) \quad \textrm{(since $f$ is linear)}\\
                                                    &= g^{-1} \xi - g^{-1} f(v) g \\
                                                    &= g^{-1} \xi - g^{-1} \rho_{J(v)} g .
\end{align*}
\end{proof}
\subsection{A Rolling Sphere on a Surface}

Now more generally, we consider a sphere rolling on a surface in $\R^3$. We will construct a connection on the trivial $SO(3)$-bundle over the surface whose parallel transport describes the rotation of the sphere. However, in this generality the connection need not be a pull-back of the natural connection on $\so(3) \times SO(3)$. Let $X$ be a smooth orientable surface in $\R^3$ and let $\n$ be the unit normal vectorfield pointing to the side on which the sphere rolls. Let $J$ be the automorphism of the tangent bundle  $TX$ that rotates each tangent space counterclockwise $\pi/2$ with axis of rotation given by the unit normal $\n$. With the natural identification of $T_xX$ with $\R^2$, $J_x(v)=\n(x) \times v$. As compared with the planar surface it is now more complicated to describe the infinitesimal rotation at a point $x \in X$ in the direction $v \in T_x X$ because of the twisting and turning of the tangent spaces of the surface. 

The unit normal vectorfield is a map $\n: X \ra \R^3$. The derivative of $\n$ at $x$ is a linear map $D\n(x): T_x X \ra \R^3$. Differentiating the constant function  $\langle \n(x),\n(x) \rangle=1$ shows that $\langle \n(x),D\n(x)(v) \rangle =0$. Therefore, $D\n(x)(v) \in T_x X$ and hence $v + D\n(x)(v)$ also lies in $T_xX$. 
Then the vector $J(v + D\n(x)(v))$ is the axis of the infinitesimal rotation of the sphere at $x$ in the direction $v$.  We define a connection on the trivial bundle $X \times SO(3)$ whose horizontal space at $(x,g)$ is
\[  H_{(x,g)}=\{(v, \rho_{J(v + D\n(x)(v))}g) | v \in T_x X\} \subset T_x X \times T_g SO(3). \]
The connection 1-form $\alpha$ is given by
\[ \alpha_{(x,g)}(v,\xi)= g^{-1}\xi- g^{-1}\rho_{J(v + D\n(x)(v))}g .\]
Let $\omega = \iota^* \alpha$ be the local connection 1-form. Thus,
\[ \omega_x(v)=\alpha_{(x,I)}(v,0)=-\rho_{J(v + D\n(x)(v))}.\]

\subsection{A Rolling Sphere on a Sphere}
When the surface $X$ is itself a sphere it is possible to explicitly compute the local curvature form. Let $X$ be the sphere of radius $r$ centered at the origin. Then $\n(x)=x/r$, $D\n(x)(v)=v/r$, and $v+D\n(x)(v) = (1 +1/r)v$. Recall that $J_x(v)=\n(x) \times v$, and so 
\[  J(v + D\n(x)(v))=\frac{x}{r} \times (v +\frac{v}{r}) =\frac{1}{r}\left(1+\frac{1}{r} \right) ( x \times v). \]
In order to compute the local curvature form $d \omega + \frac{1}{2} [ \omega, \omega]$ in coordinates we use the isomorphism $\rho$ between $\R^3$ and $\so(3)$ in order to treat $\omega$ as an $\R^3$-valued 1-form. Hence, 
\[ \omega_x(v)= -\frac{1}{r}\left(1+\frac{1}{r} \right) ( x \times v).\]
With coordinates $x=(x_1,x_2,x_3)$ and $v=(v_1,v_2,v_3)$, 
\[ \omega_x(v)=-\frac{1}{r}\left(1+\frac{1}{r} \right)(x_2v_3-x_3 v_2,x_3 v_1-x_1 v_3, x_1 v_2-x_2 v_1) .\]
Hence, 
\begin{align*} \omega &= -\frac{1}{r}\left(1+\frac{1}{r} \right)(x_2dx_3-x_3 dx_2,x_3 dx_1-x_1 dx_3, x_1 dx_2-x_2 dx_1) \\
  d\omega&= -\frac{2}{r}\left(1+\frac{1}{r} \right)(dx_2 \wedge dx_3,dx_3 \wedge dx_1,dx_1 \wedge dx_2) .
\end{align*}
Evaluating $d\omega$ on a pair of tangent vectors $u, v$ gives
\[ d \omega (u,v)= -\frac{2}{r}\left(1+\frac{1}{r} \right) u \times v .\]
Next we evaluate $\frac{1}{2} [\omega,\omega]$. Recall that $[\omega,\omega]$ is defined so that its value on a pair of tangent vectors $u$ and $v$ is
\[ [\omega,\omega](u,v) = [\omega(u),\omega(v)]-[\omega(v),\omega(u)] =2 [\omega(u),\omega(v)].\]
Therefore, $\frac{1}{2} [\omega,\omega](u,v)=[\omega(u),\omega(v)]$.
In this case the bracket operation is the cross product and so
\begin{align*} \frac{1}{2} [\omega,\omega]_x(u,v) &= \omega_x(u) \times \omega_x(v) \\
   &= \left(1 + \frac{1}{r}\right)^2 \left(\frac{x}{r} \times u\right) \times \left(\frac{x}{r} \times v\right) \\
   &= \left(1 + \frac{1}{r}\right)^2 (u \times v)
\end{align*}
For the last step above we use the geometry of the cross product with $x/r$ being a unit vector normal to both $u$ and $v$ to conclude that $(x/r \times u)\times(x/r \times v) = u \times v$.
Putting the pieces together we see that 
\begin{align*} (d\omega + \frac{1}{2}[\omega,\omega] )(u,v)
=\left(1 - \frac{1}{r^2} \right) u \times v
\end{align*}

As the radius goes to infinity the curvature form approaches the curvature form for the sphere rolling on a plane--as one would expect. For the sphere rolling on the outside of a sphere of radius 1, the curvature vanishes and so the connection is flat and the holonomy is trivial around any null-homotopic path and hence around any closed path because $S^2$ is simply-connected. Fix a basepoint $x_0 \in S^2$. There is a global section of the trivial bundle $S^2 \times SO(3) \ra S^2$ mapping $x \in S^2$ to the holonomy along any path from $x_0$ to $x$. Antipodal points take the same value as can be seen by rolling the sphere along a great circle from the north to south pole. This section is an integral surface for the horizontal distribution of the connection. It is possible to describe this map in coordinates explicitly and it is especially nice using unit quaternions, i.e., $S^3$, to represent rotations. The quaternion $q$ defines a rotation by mapping $v=(v_1,v_2,v_3)\in \R^3$ to $q(iv_1+jv_2+kv_3)q^{-1}$. The Lie group homomorphism $S^3 \ra SO(3)$ has kernel $\{\pm 1\}$. In fact, $S^3$ is the universal cover of $SO(3)$ and this homomorphism is the projection. Using $(0,0,1)$ for the basepoint, the holonomy map from $S^2$ to $SO(3)$ lifts to a map to $S^3$ given by $(x,y,z) \ra z -iy +xj$. This represents the counter-clockwise rotation through the angle $2 \cos^{-1} z$ about the axis $(-y,x,0)$. One can easily check that rolling the unit sphere from the north pole $(0,0,1)$ to $(x,y,z)$ does indeed turn the sphere through twice the angle between the two points and with axis that is normal to $(x,y,0)$.

To roll the sphere on the inside of a sphere of radius $r$ simply change the unit normal to $-x/r$ and follow the same calculations. The local curvature form $d\omega + \frac{1}{2}[\omega,\omega]$ turns out to be exactly the same
\begin{align*} (d\omega + \frac{1}{2}[\omega,\omega] )(u,v)
=\left(1 - \frac{1}{r^2} \right) u \times v .
\end{align*}
Although the curvature forms are equal, the connections are not the same and parallel transport is not the same. This can be seen easily for $r=1$. Rolling inside the sphere produces no movement at all; the rolling sphere stays fixed. Rolling the sphere on the outside does change the configuration along non-constant paths.

\textbf{Acknowledgment.} 
Thanks to Jim Delany for helpful discussions about the rolling sphere and for Mathematica code to compute holonomy for the plane and spheres.


\end{document}